\pdfoutput=1
\documentclass[microtype]{gtpart}
\usepackage[margin=1.35in]{geometry}

\newcommand{\br}{\mathbb{R}}

\newcommand{\bz}{\mathbb Z}
\newcommand{\bn}{\mathbb N}

\newcommand{\ssm}{\smallsetminus}

\DeclareMathOperator{\mcg}{MCG}
\DeclareMathOperator{\pmcg}{PMCG}

\DeclareMathOperator{\Homeo}{Homeo}

\newtheorem{Thm}{Theorem}
\newtheorem{Thm*}{Theorem}

\newtheorem{Cor}[Thm]{Corollary}
\newtheorem{Cor*}[Thm*]{Corollary}

\theoremstyle{definition}

\numberwithin{equation}{section}

\title{Three perfect mapping class groups}

\author{Nicholas G. Vlamis}
\address{Department of Mathematics \\ CUNY Queens College \\ Flushing, NY 11367}
\email{nicholas.vlamis@qc.cuny.edu}

\begin{document}  

\begin{abstract}
We prove that the mapping class group of a surface obtained from removing a Cantor set from either the 2-sphere, the plane, or the interior of the closed 2-disk has no proper countable-index normal subgroups.  
The proof is an application of the automatic continuity of these groups, which was established by Mann.
As corollaries, we see that these groups do not contain any proper finite-index subgroups and that each of these groups have trivial abelianization.
\end{abstract}

\maketitle

\vspace{-0.39in}


\section{Introduction}

The \emph{mapping class group}, denoted \( \mcg(S) \), of a topological surface \( S \) is the group of isotopy classes of orientation-preserving homeomorphisms \( S \to S \) where all homeomorphisms and isotopies are required to fix the boundary of \( S \) pointwise.
We can therefore view \( \mcg(S) \) as a quotient of a subgroup of the homeomorphism group of \( S \).
The homeomorphism group of \( S \) equipped with the compact-open topology is a topological group; this allows us to view \( \mcg(S) \) as a topological group by giving it the associated quotient topology (see \cite{AramayonaBig} for a discussion of the topological properties of mapping class groups). 

A surface is of \emph{finite type} if it is homeomorphic to a compact surface with finitely many points removed; otherwise, it is of \emph{infinite type}. 
It is known that the mapping class group of a finite-type surface is residually finite \cite{GrossmanResidual} and that every finite group appears as a quotient of some finite-index subgroup of the mapping class group of each closed surface of positive genus \cite{MasbaumAll}.
In particular, mapping class groups of finite-type surfaces are rich in finite-index subgroups.

In contrast, the subgroup structure of mapping class groups of infinite-type surfaces is not well understood and relatively few finite-index subgroups are known.
In general, finite-index subgroups are sparse as mapping class groups of infinite-type surfaces fail to be residually finite \cite{PatelAlgebraic}.
The main examples of finite-index subgroups the author is aware of come from the action of the mapping class group on the space of topological ends of the underlying surface or from pulling back finite-index subgroups of mapping class groups of closed surfaces under a homomorphism induced by embedding an infinite-type (finite-genus) surface into a closed surface.
But there are many surfaces where these constructions produce no finite- or countable-index subgroups.

In this short note, we exhibit three mapping class groups without any  proper countable-index normal subgroups.

\begin{Thm}
\label{thm:main}
The mapping class group of a surface obtained by removing a Cantor set from either the 2-sphere, the plane, or the interior of the closed 2-disk does not contain a proper countable-index normal subgroup.
\end{Thm}

A topological group \( G \) has the \emph{automatic continuity property} if every abstract homomorphism from \( G \) to a separable topological group is continuous.
The proof of Theorem~\ref{thm:main} is an application of a result of Mann \cite{MannAutomatic} showing that the mapping class groups of the 2-sphere minus a Cantor set and the plane minus a Cantor set have the automatic continuity property\footnote{There is a slight discrepancy in that Mann allows for mapping classes to be orientation-reversing.  However, she proves automatic continuity by showing that these groups have a stronger property (they are Steinhaus) that is inherited by open, normal, countable-index subgroups and hence her results apply to the orientation-preserving subgroup.}.
There are relatively few known mapping class groups with the automatic continuity property---and it is an open problem to classify those that do---and it is this lack of knowledge that prevents the proof of Theorem~\ref{thm:main} from applying to a larger class of surfaces.

Calegari and Chen have informed the author that they have independently established Theorem~\ref{thm:main}, which is to appear in forthcoming work \cite{CalegariNormal}.
The proof strategy of Calegari and Chen is fundamentally different and does not use automatic continuity.

\subsection*{Corollaries}

If \( H \) is a proper finite-index subgroup of a group \( G \), then \( G \) acts on the left cosets of \( H \), which induces a nontrivial homomorphism from \( G \) to  a finite symmetric group; hence, the kernel of this action is finite index in \( G \).
Therefore, if a group contains a proper finite-index subgroup, then it contains a proper finite-index normal subgroup.
This yields the following corollary: 

\begin{Cor}
The mapping class group of a surface obtained by removing a Cantor set from either the 2-sphere, the plane, or the interior of the closed 2-disk does not contain a proper finite-index subgroup.
\end{Cor}

Note, however, that any nontrivial mapping class group has proper countable-index subgroups:
the mapping class group acts on a number of countable sets---e.g. the isotopy classes of closed curves and the surface's first homology group with integer coefficients---and the stabilizer of any element under such an action has countable index. 

Our second corollary involves abelian quotients.
Every abelian group can be embedded into a direct sum of countable abelian groups: this follows from that fact that every abelian group can be embedded as a subgroup of a divisible group (see \cite[Theorem 24.1]{Fuchs}) and that every divisible group is a direct sum of copies of the rationals and copies of quasicyclic groups (see \cite[Theorem 23.1]{Fuchs}).
Hence, taking the projection to a factor that intersects the image of the embedding nontrivially, we can conclude that every nontrivial abelian group has a nontrivial countable quotient.
Therefore, a group in which every countable quotient is trivial must be perfect---that is, it must have trivial abelianization, or equivalently, every abelian quotient is trivial.
Together with Theorem~\ref{thm:main}, this proves:

\begin{Cor}
\label{cor:perfect}
The mapping class group of a surface obtained by removing a Cantor set from either the 2-sphere, the plane, or the interior of the closed 2-disk is perfect.
\end{Cor}

Calegari \cite{CalegariBig2} has previously shown that the group of orientation-preserving homeomorphisms and the mapping class group of the 2-sphere with a Cantor set removed are perfect; in fact, he showed that these groups are uniformly perfect---more precisely, Calegari showed that every homeomorphism can be expressed as a product of less than four commutators.  
The fact that Calegari's proof fails for the plane minus a Cantor set was the impetus for the work of Bavard~\cite{BavardHyperbolic} (see \cite{BavardHyperbolic2} for an English translation).
In particular, Bavard constructs quasimorphisms of the mapping class group of the plane minus a Cantor set, which shows that the group cannot be uniformly perfect---answering a question posed by Calegari.
Bavard's work plays a large role in the recent interest in mapping class groups of infinite-type surfaces (see \cite{AramayonaBig} for a survey). 

As a final corollary, we can state the above results in the setting of homeomorphism groups.
Mann \cite{MannAutomatic} also showed that the group of orientation-preserving homeomorphisms of the sphere minus a Cantor set and of the plane minus the Cantor set have the automatic continuity property.
Therefore, given an epimorphism from one of these groups to a countable group, we can view this epimorphism as a continuous homomorphism by equipping the codomain with the discrete topology and using the automatic continuity property of the domain. 
The continuity of the homomorphism together with the discreteness of the codomain imply that any two isotopic homeomorphisms must have the same image; hence, the given epimorphism factors through the mapping class group.
We can now apply Theorem~\ref{thm:main} and its subsequent corollaries to obtain:

\begin{Cor}
The group of orientation-preserving homeomorphisms of a surface obtained by removing a Cantor set from  either the 2-sphere or the plane does not contain a proper countable-index normal subgroup; hence, both of these groups are perfect and neither contains a proper finite-index subgroup.
\end{Cor}

\subsection*{Motivation and context}
The pure mapping class group of a surface \( S \), denoted \( \pmcg(S) \),  is the subgroup of  \( \mcg(S) \) obtained as the kernel of the action of \( \mcg(S) \) on the space of topological ends of \( S \). 
Together with Aramayona and Patel, the author in \cite{AramayonaFirst} showed the existence of nontrivial homomorphisms \( \pmcg(S) \to \bz \) whenever \( S \) has at least two nonplanar ends.
In the original version of that article, it was conjectured that all other pure mapping class groups of surfaces of genus at least two were perfect (see \cite{DomatFirst} for a discussion on genus-1 and genus-0 surfaces).
The motivation of that conjecture was a result of Patel and the author \cite{PatelAlgebraic} showing that in those cases the pure mapping class group had a dense perfect subgroup (the perfectness follows from the fact that pure mapping class groups of finite-type surfaces of genus at least three are perfect, see \cite[Theorem~5.2]{Primer}).
Domat (with an appendix with Dickmann, which is relevant in the case of the Loch Ness monster below)~\cite{DomatBig} disproved this conjecture by showing the existence of homomorphisms from pure mapping class groups to the rationals whenever the underlying surface is of infinite type and has at most one nonplanar end.

Now, there are two infinite-type surfaces without boundary in which the mapping class group is equal to the pure mapping class group, namely the Loch Ness monster surface---that is, the borderless one-ended infinite-genus surface---and the surface obtained from removing a single point from the Loch Ness monster surface.  
In both cases, Domat and Dickmann's result holds and there exists an epimorphsm from each of their mapping class groups to the rationals.
Again, since this group contains a dense perfect subgroup, Domat and Dickmann's homomorphisms necessarily fail to be continuous.
The construction of such a homomorphism must---directly or indirectly---invoke  the Axiom of Choice (this is a general fact about homomorphisms of Polish groups, see \cite{RosendalAutomatic} for a discussion), so they are naturally a bit mysterious.
However, Mann \cite{MannAutomatic} has proven that every homomorphism from any of the groups appearing in Theorem~\ref{thm:main} to a separable topological group is continuous. 
In particular, the type of homomorphism constructed by Domat and Dickmann cannot exist in these cases, and sure enough, we see that these groups are in fact perfect (Corollary~\ref{cor:perfect}).

We would also like to note that the mapping class groups of the sphere minus a Cantor set and the plane minus the Cantor set  differ drastically from a geometric viewpoint.  
First, we note that Mann and Rafi \cite{MannLarge} have shown that both the mapping class group of the sphere minus a Cantor set and the plane minus a Cantor set have a canonical metric that is well-defined up to quasi-isometry.
In this metric, Mann and Rafi have shown that the mapping class group of the sphere minus a Cantor set has bounded diameter (and in fact, every continuous left-invariant pseudo-metric on the group has bounded diameter).
In contrast, as follows from Bavard's work, the mapping class group of the plane minus a Cantor set admits infinite-diameter continuous left-invariant pseudo-metrics and hence is infinite-diameter in its canonical metric.
Moreover, Schaffer-Cohen has shown that the mapping class group of a plane minus a Cantor set equipped with its canonical metric is a Gromov hyperbolic metric space \cite{SchafferGraphs}. 

Given this discussion, from the geometric picture, it is natural to suspect that the mapping class group of the sphere minus a Cantor set should not have any quotients with interesting geometry, such as infinite finitely generated groups; however, there is no such reason to suspect this to be true of the mapping class group of the plane minus a Cantor set.
Nonetheless, this is a consequence of Theorem~\ref{thm:main}.

\subsection*{Acknowledgments}

The author is grateful to Danny Calegari, Tyrone Ghaswala, and Jes\'us Hern\'andez Hern\'andez for their comments on an earlier version.
The author recognizes support from PSC-CUNY Award \#63524-00 51.


\section{Proof of Theorem~\ref{thm:main}}

Let \( \mathbb D^2 \) denote the closed 2-disk and \( \mathbb S^2 \) the 2-sphere.
Let us view \( \mathbb D^2 \subset \br^2 \subset \mathbb S^2 \) and let \( C \) denote an embedded copy of the Cantor set in the interior of \( \mathbb D^2 \).
Let \( \Sigma \) denote either \( \mathbb S^2 \ssm C \) or \( \br^2\ssm C\) and let \( H \) be a countable group such that there exists an epimorphism \( \Phi \co \mcg(\Sigma) \to H \).
We will show that \( H \) must be the trivial group.
(The case of the disk will be handled separately at the end.)

Let \( \{c_n\}_{n\in\bn} \) denote a sequence of nontrivial isotopy classes of simple closed curves on \( \Sigma \) with the following property: for every compact subset \( K \) of \( \Sigma \) there exists \( N_K \in \bn \) such that \( c_n \) has a representative disjoint from \( K \) for every \( n > N_K \).
Let \( T_n \in \mcg(\Sigma) \) denote the left Dehn twist about \( c_n \).
By the choice of \( \{c_n \}_{n\in\bn} \), we have that \( \{T_n\}_{n\in\bn} \) is a convergent sequence in \( \mcg(\Sigma) \) limiting to the identity (this readily follows from the definition of the compact-open topology and of a Dehn twist). 
Now, with the discrete topology, \( H \) is a separable topological group, and by Mann \cite{MannAutomatic}, every homomorphism from \( \mcg(\Sigma) \) to a separable topological group is continuous; hence, \( \Phi \) is continuous.
In particular, there exists \( M \in \bn \) such that \( \Phi(T_n) \) is trivial for all \( n > M \).

Together with the universality of the Cantor set, the classification of (infinite-type) surfaces (see \cite{Richards}) implies that given any two simple closed curves on \( \Sigma \), say \( \gamma \) and \( \delta \), neither of which bound a disk nor a once-punctured disk, then  there exists a homeomorphism \( f \co \Sigma \to \Sigma \) satisfying \( f(\gamma) = \delta \) (this is the so-called \emph{change of coordinates principle} in \cite[Section 1.3]{Primer}).
It follows that any two left Dehn twists in \( \mcg(\Sigma) \) are conjugate and hence \( \Phi(T) \) is trivial for every Dehn twist \( T \in \mcg(\Sigma) \). 

By Patel--Vlamis \cite{PatelAlgebraic}, the closure of the subgroup of \( \mcg(\Sigma) \) generated by Dehn twists is \( \pmcg(\Sigma) \); hence, by the continuity of \( \Phi \), the kernel of \( \Phi \) contains \( \pmcg(\Sigma) \). 
It can be deduced from Richards's proof of the classification of surfaces in \cite{Richards} that \( \mcg(\Sigma)/\pmcg(\Sigma) \) is isomorphic to the homeomorphism group of the space of ends of \( \Sigma \) (this was pointed out to the author by S.~Afton), which in this case is \( \Homeo(C) \), the group of homeomorphisms the Cantor set.
Now, the group \( \Homeo(C) \) is simple \cite{AndersonAlgebraic} and hence \( \mcg(\Sigma)/\pmcg(\Sigma) \) is a simple group. 
Therefore, as \( \pmcg(\Sigma) < \ker \Phi \), if \( H \) is nontrivial, then \( \Phi \) factors through \( \mcg(\Sigma)/\pmcg(\Sigma) \) and \( H \) must be  isomorphic to \( \Homeo(C) \); but, then \( H \) is uncountable, which contradicts the countability assumption on \( H \).
Thus, \( H \) is trivial as desired.

Now, assume that \( \Phi \co \mcg(\mathbb D^2 \ssm C) \to H \) is an epimorphism and \( H \) is countable.
Let \( Z \) denote the infinite cyclic group generated by the left Dehn twist \( T \) about the boundary component of \( \mathbb D^2 \ssm C \).  
The embedding \( \mathbb D^2 \ssm C \hookrightarrow \br^2 \ssm C \) induces a homomorphism from \( \mcg(\mathbb D^2 \ssm C) \) to \( \mcg(\mathbb R^2 \ssm C) \) whose kernel is \( Z \).
The proof of this fact is identical to the finite-type case given in \cite[Theorem 3.18]{Primer} as long as one notes that the Alexander method holds for infinite-type surfaces \cite{HernandezAlexander}.
This allows us to see that \( \Phi \) induces an epimorphism \( \mcg(\br^2 \ssm C) \to H/\Phi(Z) \).
But, every such homomorphism is trivial; hence, \( H  = \Phi(Z) \) and is therefore cyclic. 

Using the lantern relation (see \cite[Section 5.1.1]{Primer}), we can write \( T = T_1T_2T_3(T_4T_5T_6)^{-1} \), where \( T_i \) is a left Dehn twist about a non-peripheral simple closed curve for each \( i \in \{1, \ldots, 6 \} \).
But any two left Dehn twists about  non-peripheral simple closed curves in \( \mcg(\mathbb D^2 \ssm C) \) are conjugate (for the same reason as in the earlier cases), and  since \( H \) is abelian, any two left Dehn twists have the same image under \( \Phi \).
Therefore, we see that \( \Phi(T) \), and hence \( \Phi(Z) \), is trivial; in particular, \( H \) is trivial.
\qed

\bibliographystyle{plain}
\bibliography{references}

\end{document}